\documentclass[12pt]{amsart}

\usepackage{amsmath}
\usepackage{amssymb}
\usepackage{latexsym}
\newtheorem{theorem}{Theorem}[section]

\newtheorem{prop}[theorem]{Proposition}

\newenvironment{proof*}{\vskip 2mm\noindent {}}{\hfill $\Box$ \vskip 2mm}
\numberwithin{equation}{section}
\newcommand{\C}{{\mathbb{C}}}
\newcommand{\B}{{\mathbb{B}}}

\newcommand{\D}{{\mathbb{D}}}

\newcommand{\eps}{\varepsilon}

\begin{document}

\title{Green vs. Lempert functions: a minimal example}
                     
\author{Pascal J. Thomas}

\keywords{pluricomplex Green function, Lempert function, analytic disks, Schwarz Lemma}
\subjclass{[2010]32U35, 32F45}

\address{Pascal J. Thomas\\ 
Universit\'e de Toulouse\\ UPS, INSA, UT1, UTM \\
Institut de Math\'e\-matiques de Toulouse\\
F-31062 Toulouse, France} \email{pthomas@math.univ-toulouse.fr}

\begin{abstract}
The Lempert function for a set of poles in a domain of $\mathbb C^n$ at a point $z$ is obtained by
taking a certain infimum over all analytic disks going through the poles and the point $z$,
and majorizes the corresponding multi-pole pluricomplex Green function. Coman proved that both coincide in the case of sets of two poles in
the unit ball.  We give an example of a set of three poles in the unit ball where this equality fails. 
\end{abstract}

\thanks{Part of this work was done when the author was a guest of the Hanoi University of Education in March 2011. He wishes to thank his colleagues there for their hospitality, in particular Do Duc Thai and Nguyen Van Trao. A detailed exposition of the results in a workshop organized in Bedlewo by colleagues from the Jagiellonian University in Cracow resulted in some streamlining of the proof. Finally, the Semester in Complex Analysis and Spectral Theory of the Centre de Recerca Matem\`atica at the Universitat Aut\`onoma de Barcelona provided the occasion to mention the result in a talk, and to put the finishing touches on the submitted manuscript.}

\maketitle

\section{Introduction}
Let $\Omega $ be a domain in $\C^n$, and $a_j \in \Omega$, $j=1,...,N$.
The pluricomplex Green function with logarithmic singularities at 
$S:=\{a_1,\dots,a_N\}$ is defined by 
$$
G_S(z) := \sup  \left\lbrace u \in PSH(\Omega, \mathbb R_-) :  u(z) \le \log |z-a_j|+C_j, j=1,...,N \right\rbrace ,
$$
where $PSH(\Omega, \mathbb R_-)$ stands for the set of all negative plurisubharmonic functions in $\Omega$.  When $\Omega$ is hyperconvex, this solves the Monge-Amp\`ere equation with
right hand side equal to $\sum_{i=1}^N \delta_{a_j}$. 

Pluricomplex Green functions
have been
studied by many authors at different levels of generality.  See
e.g.~Demailly \cite{De1}, Zahariuta \cite{Za}, Lempert \cite{Lempert}, Lelong \cite{Lel}, L\'arusson and Sigurdsson \cite{La-Si}.

A deep result due to Poletsky \cite{Polet}, see also \cite{La-Si}, \cite{Edig3}, is that the Green function may be computed from analytic disks:
\begin{multline}
\label{poletthm}
G_S (z)=
\inf \big\{ \sum_{\alpha: \varphi(\alpha) \in S} \log|\alpha|:
  \mbox{  such that there exists } 
  \\
   \varphi\in \mathcal {O}(\mathbb D,\Omega), \varphi(0)=z
\big\} .
\end{multline}

However, it is
tempting to pick only one $\alpha_j \in \varphi^{-1}(a_j)$, $1\le j \le N$, 
which motivated the definition of
Coman's Lempert function 
   \cite{Coman}:
\begin{multline}
\ell_S (z):=
\inf \big\{ \sum^N_{j=1} \log|\zeta_j|:
  \varphi(0)=z,
\\
\varphi(\zeta_j)=a_j, j=1,...,N
\mbox{  for some } \varphi\in \mathcal {O}(\mathbb D,\Omega) \big\} ,
\end{multline}
where $\mathbb D$ is the unit disc in $\mathbb C$. 

One easily sees that $\ell_S (z) \ge G_S(z)$ without recourse to 
\eqref{poletthm}; the fact 
that equality holds when $N=1$ and $\Omega$ is convex is part of Lempert's 
celebrated theorem, which was, in fact, the starting point for
many of the notions defined above \cite{Lempert},
see also \cite{Edig2}. Coman proved that equality holds when $N=2$ and
$\Omega=\mathbb B^2$, the unit ball of $\mathbb C^2$ \cite{Coman}.  The goal of this note is to present an example that shows that this is as far as it can go.

\begin{theorem}
\label{main}
There exists a set of $3$ points $S \subset \mathbb B^2$ such that
for some $z \in \mathbb B^2$, $\ell_S (z) > G_S(z)$.
\end{theorem}

Other examples in the same direction have been found in  \cite{CarlWieg}, \cite{TraoTh}, \cite{NikoZwo}. The interesting features of this one are that it involves no multiplicities
and is minimal in the ball. Examples with an arbitrary number of points can be deduced from it.  Let $z_0\in \B^2$ satisfy $\ell_S(z_0)-G_S(z_0)=:\eps_0>0$. Consider $S':= S \cup \{a_4, \dots, a_N\}$  with all the $a_j$ close
enough to the boundary so that $\ell_{S'}(z_0) \ge \ell_S(z_0)-\eps_0/2$ 
(the Schwarz lemma shows that $|\zeta_j|\to1$ when $\varphi(\zeta_j)=a_j$ and
$|a_j|\to1$). Then $\ell_{S'}(z_0) > G_S(z_0) \ge G_{S'}(z_0)$, q.e.d. (I thank 
Nikolai Nikolov for sharing this observation with me). 

Furthermore,  the corresponding Green function can be recovered, up to a bounded error, by using an analytic disk with just one more pre-image than the number of points:
one of the points has exactly two pre-images and each of the other two points,
only one, see \cite[\S 6.8.2, Lemma 6.16]{MRST}.

More specifically, the Theorem will follow from a precise calculation in the bidisk $\D^2$.
Let $S_\eps = \left\{ (0,0), (\rho(\eps), 0), (0, \eps)\right\} \subset \D^2$,
where  $\lim_{\eps\to0} \rho(\eps)/\eps=0$.

\begin{prop}
\label{bidisk}
There exist $C_1 >0$ and, for any $\delta \in (0,\frac14)$, 
$r_0=r_0(\delta)>0$
such that 
for any $z=(z_1,z_2)\in \D^2$ with 
\begin{equation}
\label{sector}
\frac{1}{2} |z_2|^{3/2} \le |z_1| \le |z_2|^{3/2},
\quad \|z\| <r_0,
 \end{equation}
 there exists $\varepsilon_0=\varepsilon_0(z,\delta) >0$ 
such that for any $\eps$ with $|\eps|<\varepsilon_0$,
then 
\begin{equation}
\label{grup}
G_{S_\eps} (z) \le 2 \log|z_2| + C_1,
\end{equation}
\begin{equation}
\label{lembel}
\ell_{S_\eps} (z) \ge (2-\delta) \log|z_2| . 
\end{equation}
\end{prop}

\begin{proof*}{\it Proof of Theorem \ref{main}.}
If $U,V$  are domains, and $S\subset U \subset V$, then the definitions of the Green and Lempert functions imply that $G^U_S (z) \ge G^V_S (z)$,
$\ell^U_S (z) \ge \ell^V_S (z)$. For $|\eps|$ small enough, 
$S_\eps\subset \B^2$.
When $|z_1|=|z_2|^{3/2}$, so that $z$ verifies \eqref{sector},
the inclusion  $\B^2 \subset \D^2$
implies 
$$
\ell^{\B^2}_{S_\eps} (z) \ge \ell^{\D^2}_{S_\eps} (z) \ge (2-\delta) \log|z_2| .
$$
Using the fact that  $\frac{\sqrt2}{2}\D^2 \subset \B^2$ and the invariance of
the Green function under biholomorphic mappings,
$$
G^{\B^2}_{S_\eps} (z) \le G^{\frac{\sqrt2}{2}\D^2}_{S_\eps} (z) =
 G^{\D^2}_{\sqrt2 S_\eps} (\sqrt2 z) \le 
  2 \log|z_2| + \log 2 + C_1.
$$
The last inequality follows from the fact that $\sqrt2 z$ still verifies \eqref{sector},  
and $\sqrt2 S_\eps$ has the same form as $S_\eps$, so we can apply \eqref{grup}. 

Comparing the last two estimates, we see that for $|z_2|$ small enough
 and $|\eps|<\varepsilon_0$, 
$G^{\B^2}_{S_\eps} (z) < \ell^{\B^2}_{S_\eps} (z)$. 
\end{proof*}

{\bf Open Questions.}

1. This example is minimal in the ball, in terms of number of poles; what is the situation for the bidisk? Are the Green and Lempert functions equal when one takes two poles, not lying on a line parallel to the coordinate axes? Do they at least have the same order of singularity as one pole tends to the other?

2. What is the precise order of the singularity of the
limit as $\eps\to 0$ of the Lempert function in this case?  Looking at the available analytic disks that give the correct order of the singularity of the limit of the Green function, one finds $\frac32 \log|z_2|$, so one would hope that the Proposition can still be proved at least for $\delta < \frac12$. 

3. Do the analytic disks from \cite{MRST} yield the Green function itself, without any bounded error term? 

4. More generally, when one is given a finite number of points in a given bounded (hyperconvex) domain, is there a bound on the number of pre-images required to attain the Green function in the Poletsky formula?  For instance, is $4$ the largest possible number of pre-images required when looking at $3$ points in the ball?

{\bf Acknowledgements.}

I wish to thank Nguyen Van Trao for useful discussions on this topic,
and the referee for pointing out and correcting a mistake in the original exposition.

\section{Upper estimate for the Green function}

\begin{proof*}{\it Proof of Proposition \ref{bidisk} \eqref{grup}.}

The upper bound \eqref{grup} follows from \cite[\S 6.8.2, Lemma 6.16]{MRST}.  
For the reader's convenience, and since that paper is not (yet) generally available,
we repeat the proof here in the case that concerns us. 

We now construct an analytic disk 
passing twice through one of the poles. Our disk will be a perturbation of
the Neil parabola $\zeta \mapsto (\zeta^3, \zeta^2)$.

We write $s(\eps) = \rho(\eps)/\eps = o(1)$.

Choose complex numbers $\lambda, \mu$ such that
$$
\lambda^2 := \frac{z_1}{z_2(z_2-\eps)}
\left( \frac{z_1}{z_2-\eps} + s(\eps) \right) ; \quad
\mu^2 := \eps + \left( \frac{s(\eps)}{2 \lambda} \right)^2.
$$
Let
$$
\Psi_{\lambda, \mu}(\zeta) :=
\left(  \left( \lambda \zeta - \frac12 s(\eps) \right) (\zeta^2-\mu^2),
\zeta^2 -  \left( \frac{s(\eps)}{2 \lambda} \right)^2 
\right) .
$$
Then by construction
$$
\Psi_{\lambda, \mu}(\mu) = \Psi_{\lambda, \mu}(-\mu)=(0,\eps) ,
\Psi_{\lambda, \mu} ( \frac{s(\eps)}{2 \lambda}  ) = (0,0) ,
\Psi_{\lambda, \mu} ( -\frac{s(\eps)}{2 \lambda}  ) = (\eps s(\eps), 0),
$$
so we have a disk passing through all three poles of $G_{\eps}$.
Furthermore, choosing 
$$
\zeta_z := \frac1\lambda \left( \frac{z_1}{z_2-\eps} + \frac{s(\eps)}2 \right),
$$
we have $\Psi_{\lambda, \mu}(\zeta_z)= z$. Notice that
$$
\zeta_z^2 = \frac{z_2(z_2-\eps)}{z_1} \left( \frac{z_1}{z_2-\eps} + \frac{s(\eps)}2  \right)^2 \left( \frac{z_1}{z_2-\eps} + s(\eps)  \right)^{-1},
$$
so for any $\eta>0$ 
there exists $\eps_0(\delta, \eta)>0$ such that for $|\eps|< \eps_0(\delta, \eta)$,
for any $z$ such that $\delta\le \frac12 |z_2|^{3/2} \le |z_1| \le |z_2|^{3/2}\le 1$, 
\begin{equation}
\label{modzeta}
\left| |\zeta_z| - |z_2|^{1/2}\right| \le \eta. 
\end{equation}
In particular, by choosing $\eta$
small enough we ensure that $\zeta_z \in \D$.  We need a more general fact.

Claim.  

Let $\eta>0$, and $\delta >0$.  Then there exists 
 $\eps_1 = \eps_1(\delta, \eta)>0$ such that for any $\eps$
with $|\eps|\le \eps_1$, for any
 $z$ such that 
$\delta\le \frac12 |z_2|^{3/2} \le |z_1| \le |z_2|^{3/2}\le 1$, we have $\Psi_{\lambda, \mu} (D(0,1-\eta) )
\subset \D^2$. 

Proof of the Claim.

For $|\eps|\le \delta^{2/3}/2$, $|z_2|/2 \le |z_2-\eps| \le 2|z_2|$, so
$$
|\lambda|^2 \ge \left| \frac{z_1}{2z_2^2} \right| 
\left(  \left| \frac{z_1}{2z_2} \right| - |s(\eps)| \right) 
\ge \left| \frac{z_1^2}{8z_2^3} \right| \ge \frac1{32},
$$
for $\eps$ small enough. So when $|\zeta|\le 1-\eta$, 
$$
\left| \Psi_{\lambda, \mu,2}(\zeta)\right| 
\le (1-\eta)^2 + 256 |s(\eps)|^2 <1
$$
for $\eps$ small enough. 

In a similar way, given $\eta'$, 
for $\eps$ small enough depending on $\delta$ and $\eta'$, we have
$|z_2| \le (1+\eta') |z_2-\eps|$, so
$$
|\lambda|^2 \le (1+\eta')^2 \left| \frac{z_1}{z_2^2} \right| 
\left(  \left| \frac{z_1}{z_2} \right| +  \frac{|s(\eps)|}{(1+\eta') } \right) 
\le (1+\eta')^3 \left| \frac{z_1^2}{z_2^3} \right| \le (1+\eta')^3
$$
for  $\eps$ small enough. Choose $\eta'$ so that 
$(1+\eta')^3 = (1+\eta)$.
When $|\zeta|\le 1-\eta$, 
$$
\left| \Psi_{\lambda, \mu,1}(\zeta)\right| 
\le
\left(  (1+\eta)  (1-\eta) + \frac12 |s(\eps)| \right) 
\left( (1-\eta)^2 + |\eps| +  64^2 |s(\eps)|^2  \right) <1
$$
for  $\eps$ small enough.  \hfill \qed

So now the function $v(\zeta):= G_{\eps} \left( \Psi_{\lambda, \mu}((1-\eta)\zeta) \right) $
is negative and subharmonic on $\D$.  Furthermore, it has logarithmic poles at the points
$\pm \frac\mu{1-\eta}$ and $\pm  \frac{s(\eps)}{2 \lambda(1-\eta)} $ ; in the cases when $\mu=0$ or
$s(\eps)=0$, we get a double logarithmic pole at the corresponding point. 

Denote $d_{G}(\zeta, \xi):= \left| \frac{\zeta-\xi}{1-\zeta \bar \xi} \right|$
the invariant (pseudohyperbolic) distance between points of the unit disk. Then
\begin{multline*}
G_{\eps} (z) = v(\zeta_z) \le
\log d_G(\zeta_z, \frac\mu{1-\eta}) + \log d_G(\zeta_z, -\frac\mu{1-\eta}) \\
+ \log d_G(\zeta_z, \frac{s(\eps)}{2 \lambda (1-\eta)})  + 
\log d_G(\zeta_z, -\frac{s(\eps)}{2 \lambda (1-\eta)}). 
\end{multline*}
By \eqref{modzeta}, choosing
$m(\delta, \eta)$ accordingly, we have, for $|\eps|\le m$,
$G_{\eps} (z) \le 4 \log |z_2|^{1/2} +O(\eta)$. 
\end{proof*}

\section{Lower estimate for the Lempert function}

\begin{proof*}{\it Proof of Proposition \ref{bidisk} \eqref{lembel}.}
The proof of \eqref{lembel}
will follow the methods and notations of \cite{Th}.
We will make repeated use of the involutive automorphisms of the unit disk given by $\phi_a(\zeta):= \frac{a-\zeta}{1-\bar a \zeta}$ for $a\in\D$,
which exchange $0$ and $a$.
Notice that the invariant (pseudohyperbolic) distance verifies
$$
d_G (a,b) := \left| \phi_a (b) \right| =  \left| \phi_b (a) \right| .
$$

Write $\rho(\eps)= \eps s(\eps)$ with $\lim_{\eps\to0} s(\eps)=0$.  

We will assume that the conclusion fails.  
That is, for any $\delta \in (0,\frac14)$,
 there exist arbitrarily small values of $|z_2|=\max(|z_1|,|z_2|)$, 
 and $|\eps|$ such that
\begin{equation}
\label{contrhyp}
\ell_{S_\eps}(z) < (2-\delta)\log |z_2| .
\end{equation}
After applying, for each analytic disk,
an automorphism of the disk which exchanges the pre-image of $(0,0)$
and $0$,
the assumption
implies that there exists a
holomorphic map $\varphi$ from $\D$ to $\D^2$ and points $\zeta_j \in \D,$ 
depending on $z$ and $\eps$, satisfying the conditions 
\begin{equation}
\label{prob}
\left\lbrace
\begin{array}{rcl}
\varphi(0)&=&(0,0)\\
\varphi(\zeta_1)&=&(\eps s(\eps),0)\\
\varphi(\zeta_2)&=&(0,\eps)\\
\varphi(\zeta_0)&=&(z_1,z_2)
\end{array}
\right.
\end{equation}
 with
\begin{equation}
\label{hyplemp}
\log |\zeta_0| + \log |\phi_{\zeta_0}(\zeta_1)| + \log |\phi_{\zeta_0}(\zeta_2)|  \le  (2-\delta)\log |z_2| .
\end{equation}

The interpolation conditions in \eqref{prob} are equivalent to the existence of two holomorphic functions $h_1$, $h_2$ from
$\D$ to itself such that
\begin{equation*}
\varphi (\zeta) = \left( \zeta \phi_{\zeta_2} (\zeta) h_1(\zeta), \zeta \phi_{\zeta_1} (\zeta) h_2 (\zeta) \right) ,
\end{equation*}
such that furthermore
\begin{eqnarray}
\label{h11}
h_1(\zeta_1)&=& \frac{\eps s(\eps)}{\zeta_1 \phi_{\zeta_2} (\zeta_1)}=:w_1, \\
\label{h12}
h_1(\zeta_0)&=& \frac{z_1}{\zeta_0 \phi_{\zeta_2} (\zeta_0)}=:w_2, \\
\label{h21}
h_2(\zeta_2)&=& \frac{\eps}{\zeta_2 \phi_{\zeta_1} (\zeta_2)}=:w_4, \\
\label{h22}
h_2(\zeta_0)&=& \frac{z_2}{\zeta_0\phi_{\zeta_1} (\zeta_0)}=:w_3.
\end{eqnarray}

By the invariant Schwarz Lemma, the existence of a holomorphic  function $h_1$ mapping $\D$ to itself and satisfying
\eqref{h11} and \eqref{h12} is equivalent to
\begin{equation}
\label{ineqex1}
\left| w_1 \right| < 1,
\left| w_2 \right| < 1,
\mbox{ and }
d_G \left( w_1, w_2 \right) < d_G \left( \zeta_1, \zeta_0 \right) = \left| \phi_{\zeta_1} (\zeta_0) \right|.
\end{equation}
In the same way, the existence of $h_2$ is equivalent to
\begin{equation}
\label{ineqex2}
\left| w_3 \right| < 1,
\left| w_4 \right| < 1,
\mbox{ and }
d_G \left( w_3, w_4 \right) < d_G \left( \zeta_2, \zeta_0 \right) = \left| \phi_{\zeta_2} (\zeta_0) \right|.
\end{equation}

As in \cite{Th}, we start by remarking that \eqref{hyplemp} can be rewritten as
\begin{multline}
\label{startrem}
-\log |w_2|-\log |w_3|=
 \log \left| \frac{\zeta_0 \phi_{\zeta_1}(\zeta_0)}{z_2} \right| +
 \log \left| \frac{\zeta_0 \phi_{\zeta_0}(\zeta_2)}{z_1} \right|
\\
\le \log |\zeta_0| + (2-\delta) \log |z_2| - \log|z_1|  - \log|z_2|
\\
\le  \log |\zeta_0| - (\frac{1}{2}+ \delta) \log |z_2| + \log 2,
\end{multline}
by \eqref{sector}.   We can rewrite this in a more symmetric fashion:
\begin{equation}
\label{symest}
\log \frac1{|w_2|} +\log \frac1{|w_3|} + \log \frac1{|\zeta_0|}
  \le (\frac{1}{2}+ \delta) \log \frac1{|z_2|} + \log 2.
\end{equation}
Since all terms 
are positive by \eqref{ineqex1}, \eqref{ineqex2}, each of the terms on
the left hand side is bounded by the right hand side. 

We will proceed as follows: we have used the contradiction hypothesis
\eqref{hyplemp} to prove that $|\zeta_0|$ and $|w_3|$ 
are relatively big.  We will prove that $|\phi_{\zeta_2} (\zeta_0)|$ has to be relatively small, which by \eqref{ineqex2} forces $|w_4|$ to be roughly as large as $|w_3|$.  This then allows us to bound $|\phi_{\zeta_1} (\zeta_2)|$ by a quantity which becomes as small as desired when $\eps$ can be made small, hence allows us to bound $|\phi_{\zeta_1} (\zeta_0)|$ by the triangle inequality. 

The final contradiction will concern $w_2=\frac{z_1}{\zeta_0 \phi_{\zeta_2} (\zeta_0)}$. On the one hand, \eqref{symest} guarantees that it is not too small; but an explicit computation of the quotient $w_1/w_4$ shows that $w_1$ must be small, and by \eqref{ineqex1} and the estimate on 
$|\phi_{\zeta_1} (\zeta_0)|$, $|w_2|$ is small as well.

We provide the details. 
From \eqref{symest}, 
\begin{equation}
\label{w3est}
 \log |w_3| 
 \ge (\frac{1}{2}+ \delta) \log |z_2| -\log 2.
\end{equation}
From \eqref{h12} and \eqref{startrem},
\begin{multline}
\label{phi02est}
\log |\phi_{\zeta_2} (\zeta_0)| = \log \left| \frac{z_1}{\zeta_0}\right|
- \log |w_2| \\
 \le \log \left| \frac{z_1}{\zeta_0}\right|
+ \log |\zeta_0| - (\frac{1}{2}+ \delta) \log |z_2| + \log 2 \le 
(1- \delta) \log |z_2| +\log 2.
\end{multline}
Since $\delta <\frac14$,
\eqref{phi02est} and \eqref{w3est} imply that
for  $|z_2| \le r_1(\delta)$, $|\phi_{\zeta_2} (\zeta_0)| <\frac12 |w_3|$, so by \eqref{ineqex2} and the triangle inequality for $d_G$, 
\begin{equation}
\label{w4est}
|w_4| \ge \frac12 |w_3|. 
\end{equation}
We now prove that both $\zeta_1$ and $\zeta_2$ must be close to $\zeta_0$ and even closer to each other. First, since \eqref{symest}
implies that $\log|\zeta_0| \ge (\frac{1}{2}+ \delta) \log |z_2| - \log 2$, by  
\eqref{phi02est},  for $|z_2| \le r_2(\delta)$, $|\phi_{\zeta_2} (\zeta_0)| \le \frac12 |\zeta_0|$.
By the triangle inequality for $d_G$, 
\begin{equation}
\label{ze2est}
\frac12 |\zeta_0|\le |\zeta_2| \le \frac32 |\zeta_0|.
\end{equation}
On the other hand, from \eqref{symest}, 
$$
\log |w_3| + \log |\zeta_0| \ge (\frac{1}{2}+ \delta) \log |z_2|-\log 2 ,
\mbox{ i.e. } |w_3 \zeta_0| \ge \frac{1}{2} |z_2|^{\delta+1/2}.
$$
Therefore, applying \eqref{w4est} and \eqref{ze2est}, 
\begin{equation}
\label{phi12est}
\left| \phi_{\zeta_1} (\zeta_2) \right| = \left| \frac{\eps}{\zeta_2 w_4}\right| \le 4 \left| \frac{\eps}{\zeta_0 w_3}\right| \le 8 |\eps| |z_2|^{-\delta-1/2}.
\end{equation}
In particular, for 
\begin{equation}
\label{epscond}
|\eps| < \frac18 |z_2|^{3/2},
\end{equation}
this implies $| \phi_{\zeta_1} (\zeta_2) | <  |z_2|^{1-\delta}$, and by the triangle inequality, 
\begin{equation}
\label{phi10est}
\left| \phi_{\zeta_1} (\zeta_0)  \right| < 
\left| \phi_{\zeta_2} (\zeta_0)  \right| + \left| \phi_{\zeta_1} (\zeta_2) \right| < 3 |z_2|^{1-\delta}.
\end{equation}

We now establish the two (contradictory) estimates for $w_2$.  On the one hand, 
\eqref{symest} implies that
\begin{equation}
\label{w2bel}
\log |w_2|  \ge (\frac{1}{2}+ \delta) \log |z_2|-\log 2 ,
\mbox{ i.e. } |w_2| \ge \frac{1}{2} |z_2|^{\delta+1/2}.
\end{equation}
On the other hand, 
$$
\left|\frac{w_1}{w_4}\right| = 
\left| \frac{\eps s(\eps)}{\zeta_1 \phi_{\zeta_2} (\zeta_1)} 
 \frac{\zeta_2 \phi_{\zeta_1} (\zeta_2)}{\eps} \right|
 =\left| s(\eps) \frac{\zeta_2}{\zeta_1} \right| .
$$
By the triangle inequality for $d_G$, when \eqref{epscond} holds,
the lower bound in \eqref{ze2est} and the corollary to \eqref{phi12est}
imply
$$ 
|\zeta_1| \ge |\zeta_2| - \left| \phi_{\zeta_1} (\zeta_2)\right|
\ge \frac12 |\zeta_0| -  |z_2|^{1-\delta} 
\ge \frac14 |\zeta_0| 
$$
for $|z_2|$ small enough, because of \eqref{symest} again.  So finally, using the upper bound in \eqref{ze2est},
$|\frac{w_1}{w_4}| \le 6 |s(\eps)|$. We choose $\varepsilon_0 < \frac18 |z_2|^{3/2}$
so that for any $\eps$ with $|\eps|\le \varepsilon_0$,
\begin{equation}
\label{epscond2}
|s(\eps)| < |z_2|^{1-\delta}.
\end{equation}
The triangle inequality for $d_G$ and \eqref{phi10est} imply that when 
$|\eps|\le \varepsilon_0$,
$$
|w_2| \le |w_1| + \left| \phi_{\zeta_1} (\zeta_0)  \right| 
\le 6 |s(\eps)| +  3 |z_2|^{1-\delta} \le  9 |z_2|^{1-\delta}.
$$
Finally, if we choose $|z_2|\le r_0(\delta)$, with
$r_0(\delta) \le \min( r_1(\delta),r_2(\delta))$, and
$9r_0(\delta)^{1-\delta} < \frac12 r_0(\delta)^{\frac12+\delta}$,
we see that for any $\eps$ with $|\eps|\le \varepsilon_0$,
this last bound contradicts \eqref{w2bel}. 
\end{proof*}

\bibliographystyle{amsplain}

\begin{thebibliography}{ABC}


\bibitem{CarlWieg}
M. Carlehed, J. Wiegerinck,
{\it Le c\^one des fonctions
plurisousharmoniques n\'egatives et une conjecture de Coman,}
Ann. Pol. Math. {\bf 80} (2003), 93--108.

\vspace{0.2cm}

\bibitem{Coman}
D. Coman, {\it The pluricomplex Green function with two poles of the
unit ball of $\mathbb C^n$,} Pacific J. Math. {\bf 194,} no 2,
257--283 (2000).

\vspace{0.2cm}

\bibitem{De1} J.-P. Demailly, {\it Mesures de Monge-Amp\`ere et mesures
pluriharmoniques}, Math. Z. 194 (1987), 519--564.
\vspace{0.2cm}

\bibitem{Edig2}
A. Edigarian, {\it A remark on the Lempert theorem,}
 Univ. Iagel. Acta Math. {\bf 32} (1995), 83--88. 

\vspace{0.2cm}

\bibitem{Edigarian}
A. Edigarian, {\it Remarks on the pluricomplex Green function,} Univ.
Iagel. Acta Math. {\bf 37} (1999), 159--164.

\vspace{0.2cm}

\bibitem{Edig3}
A. Edigarian, {\it On definitions of the pluricomplex Green function,}
 Ann. Polon. Math. {\bf 67} (1997), no. 3, 233--246. 

\vspace{0.2cm}

\bibitem{Garnet}
J. B. Garnett, {\it Bounded Analytic Functions,} Academic Press, Inc. New
York-London, 1981.

\vspace{0.2cm}

\bibitem{La-Si} F. L\'arusson, R. Sigurdsson, {\it Plurisubharmonic functions and analytic discs on manifolds}, J. Reine Angew. Math. 501 (1998), 1--39.

\vspace{0.2cm}

\bibitem{Lel} P. Lelong, {\it Fonction de Green pluricomplexe et lemmes
de Schwarz dans les espaces de Banach}, J. Math. Pures Appl. 68 (1989), 319--347.

\vspace{0.2cm}

\bibitem{Lempert}
L. Lempert, {\it La m\'etrique de Kobayashi et la repr\'esentation des
domains sur la boule,} Bull. Soc. Soc. Math. France 109 (1981),
427--474.

\vspace{0.2cm}

\bibitem{MRST}

J. I. Magn\'usson, A. Rashkovskii, R. Sigurdsson, P. J. Thomas, {\it Limits of multipole pluricomplex Green functions,}  arXiv:1103.2296.

\vspace{0.2cm}

\bibitem{NikoZwo}
N. Nikolov, W. Zwonek, {\it On the product property for the Lempert function,}
 Complex Var. Theory Appl. {\bf 50} (2005), no. 12, 939--952. 

\vspace{0.2cm}

\bibitem{Polet}
E. A. Poletsky, {\it Holomorphic currents,} Indiana Univ. Math. J.
  42 (1993), no 1, 85--144.
\vspace{0.2cm}

\bibitem{Th} P. J. Thomas, {\it An example of limit of Lempert functions}, Vietnam Journal of Mathematics 35 (2007), no. 3, 317--330.
\vspace{0.2cm}

\bibitem{TraoTh}
P. J. Thomas, N. V. Trao, {\it Pluricomplex Green and Lempert functions
for equally weighted poles.} Ark. Mat.  41 (2003),  no. 2, 381--400.
\vspace{0.2cm}

\bibitem{Za} {V.P. Zahariuta}, {\it Spaces of analytic functions and maximal
plurisubharmonic functions.} D. Sci. Dissertation, Rostov-on-Don,
1984.


\end{thebibliography}

\end{document}